\documentclass[11pt]{amsart}
\usepackage{geometry}
\usepackage{amssymb, amsfonts, amsthm, dsfont,bm}
\usepackage{mathrsfs}

\newtheorem{thm}{\textbf{Theorem}}[section]
\newtheorem{coro}{\textbf{Corollary}}[section]

\newtheorem{lemma}{\textbf{Lemma}}[section]
\newtheorem{rem}{\textbf{Remark}}[section]

\def\nn{\mathbb{N}}
\def\rr{\mathbb{R}}
\def\ee{\mathbb{E}}

\def\x{\underline{x}}
\def\y{\underline{y}}

\def\dd{\textup{d}}

\def\F{{\mathscr F}}

\def\H{\widehat{H}}

\def\s{\sigma}
\def\var{{\rm var}}
\def\Lip{{\rm Lip}}
\def\Var{{\rm Var}}

\def\E{\mathcal{E}}

\title[Concentration bounds for entropy estimators]{
Concentration bounds for entropy estimation of one-dimensional Gibbs measures}
\author[J.-R. Chazottes and C. Maldonado]{J.-R. Chazottes, C. Maldonado}
\address{CPhT, CNRS-\'{E}cole Polytechnique, 91128 Palaiseau Cedex, France}
\address{Email address: \textit{\texttt{jeanrene@cpht.polytechnique.fr}}}
\address{Email address: \textit{\texttt{maldonado@cpht.polytechnique.fr}}}
\date{} 

\begin{document}

\begin{abstract}
We obtain bounds on fluctuations of two entropy estimators for
a class of one-dimensional Gibbs measures on the full shift. They are the consequence
of a general exponential inequality for Lipschitz functions
of $n$ variables. The first estimator is based on empirical frequencies of blocks 
scaling logarithmically with the sample length.
The second one is based on the first appearance of blocks within typical samples. 
\end{abstract}

\maketitle
\tableofcontents

\section{Introduction}

Given a `sample' $x_0,x_1,\ldots$ of a finite-valued discrete-time ergodic process $\{X_n; n\in\nn\}$,
there are several ways to consistently estimate its entropy. In this paper
we shall study two estimators. One is based on empirical frequencies of blocks and is
referred to as the `plug-in' estimator. The other one is based on the first appearance or
repetition of blocks within the sample. We refer to \cite{shields} for their basic properties.
Here we are concerned with the fluctuation properties of these estimators. We will further
assume that the joint distribution of the process $\{X_n; n\in\nn\}$ is Gibbsian in a way made precise
below.

Fluctuations of the plug-in estimator were already studied in \cite{GGG} and \cite{jrdavide} from
the viewpoint of classical limit theorems. Namely, in \cite{GGG} the authors prove a central limit theorem
and in \cite{jrdavide} a large deviation principle is obtained.

Regarding the return-time and the hitting-time estimators, previous results are found in 
\cite{CGS} and \cite{ChU}. Central limit theorems and large deviations principle are
established in these papers. In the present article, we only study the hitting-time estimator.

Our aim is to obtain bounds on the fluctuations of the plug-in and hitting-time
estimator in the spirit of concentration inequalities.
Concentration inequalities became recently a widespread powerful tool in many fields of pure and applied probabilities,
as well as in functional analysis, combinatorics, computer science, etc; see for instance \cite{Led} and
\cite{Panconesi}.
In the context of dynamical systems, the first result was proved in \cite{CMS} where several applications
are presented (see also \cite{ChCS2}). Namely, an exponential inequality is proved for any separately Lipschitz function of $n$ variables for a class of piecewise expanding maps of the interval.
In our setting the same inequality holds. The proof is the same as in \cite{CMS}. It is in fact
simpler since no Markov partition is assumed therein.
In this paper we apply this exponential inequality to get some fluctuation bounds on our
entropy estimators. The only previous work where this is done for the plug-in estimator is found in \cite{AK} in
the case where the $X_i$'s are independent identically distributed random variables
taking values in a countable set. For the hitting-time estimator, no such bounds were known before, even
in the case of independent random variables.

Our main results are theorems \ref{main1}, \ref{main2} and \ref{thmWaiting}. Let us emphasize
that we establish bounds for every $n$, $n$ being the sample length, whereas the results obtained
in \cite{jrdavide,ChU,CGS,GGG} are in some sense finer but they are only asymptotic. Therefore,
our work complement the picture on the fluctuations of these entropy estimators.
Theorems \ref{main1} and  \ref{main2} concern the plug-in estimator. Theorem \ref{thmWaiting}
is about the hitting-time estimator. It should be noted that the route to get the bounds
is not as direct as for the plug-in estimator because the hitting-time a priori behaves badly.
The trick is to take advantage of its approximation by the inverse measure of the corresponding
cylinder. This is where Gibbsianness is crucial.

This paper is organized as follows. In Section \ref{setting} we recall some definitions and facts, and
state the exponential bound from which concentration inequalities follow. Section \ref{clementine}
contains our results on the plug-in estimators and the hitting-time estimator. 

\section{Setting}\label{setting}

After fixing some notations, we recall a few facts about entropy and Gibbs measures.

\subsection{Notations and definitions}

We consider the set $\Omega=A^{\nn}$ of infinite sequences $\x$ of symbols from
the finite set $A$: $\x=x_0,x_1,\ldots$ where $x_j\in A$. 
We denote by $\s$ the shift map on $\Omega$: $(\s\x)_{i} = x_{i+1}$, for all $i=0,1,\ldots$. \newline
We equip $\Omega$ with the usual distance: fix $\theta\in(0,1)$ and
for $\x\neq \y$, let $d_{\theta}(\x,\y) = \theta^{N}$ where $N$ is the largest nonnegative
integer with $x_i=y_i$ for every $0\leq i<N$. (By convention, if $\x=\y$ then $N=\infty$ and $\theta^\infty=0$, 
while if $x_0\neq y_0$ then $N=0$.)
With this distance $\Omega$ is a compact metric space. \newline
For a given string $a_0^{k-1}=a_0,\ldots,a_{k-1}$ ($a_i\in A$), the set
$[a_0^{k-1}] =\{\x\in\Omega : x_i=a_i, i=1,\ldots,k-1\}$ is the cylinder
of length $k$ based on $a_0,\ldots,a_{k-1}$. \newline
For a continuous function $f:\Omega\to\rr$ and $m\geq 0$ we define
$$
\var_{m}(f) : = \sup\{\lvert f(\x) -f(\y) \rvert :  x_{i} = y_{i}, \ i=0,\ldots,m\}\cdot
$$
It is easy to see that $\lvert f(\x) -f(\y) \rvert \leq C d_{\theta}(\x,\y)$ if and only
if $\var_{m}(f) \leq C \theta^m$, $m=0,1,\ldots$.
Let
$$
\F_\theta=\big\{f: f\;\textup{continuous},\;\var_{m}(f)\leq C\theta^m,\ m=0,1,\ldots,\ \textup{for some}\;C>0\big\}.
$$
This is the space of Lipschitz functions with respect to the distance $d_\theta$.
For $f\in \F_\theta$ let
$|f|_\theta=\sup\left\{ \frac{\var_m(f)}{\theta^m}:m\geq 0\right\}$.
We notice that $|f|_\theta$ is merely the least Lipschitz constant of $f$.

\subsection{Entropy}

Let $\nu$ be a shift-invariant probability measure on $\Omega$ and  
\[
H_{k}(\nu) = -\sum_{a_{0}^{k-1}\in A^k}\nu([a_{0}^{k-1}])\log\nu([a_{0}^{k-1}]),
\]
its `$k$-block entropy'.
Then the entropy of $\nu$ is
$$
h(\nu)=\lim_{k\to\infty} \frac{H_k(\nu)}{k}\cdot
$$
Recall that $0\leq h(\nu)\leq \log|A|$.
\subsection{Gibbs measures}

Full details for this section can be found in \cite{Bow}. Let $\phi\in\F_\theta$ and $\mu_\phi$ the associated Gibbs measure. It is the unique shift-invariant probability measure for which one can find constants $C=C(\phi)>1$
and $P=P(\phi)$ such that
\begin{equation}\label{BowenGibbs}
C^{-1}\leq
\frac{\mu_\phi\big(\big\{\y : y_{i}=x_{i}, \forall \ i\in [0,m) \big\}\big)}{\exp\left(-Pm + \sum_{k=0}^{m-1}\phi(\s^{k}\x)\right)}
\leq C
\end{equation}
for every $\x\in\Omega$ and $m\geq 1$. The constant $P$ is the topological pressure of $\phi$.
We can always assume that $P=0$ by considering the potential $\phi-P$ which yields the same
Gibbs measure.

The Gibbs measure $\mu_\phi$ satisfies the variational principle, namely 
$$
\sup\left\{h(\eta) +\int\phi \dd\eta: \eta\;\textup{shift-invariant}\right\}=h(\mu_\phi)+\int \phi \dd\mu_\phi=P=0.
$$
More precisely, $\mu_\phi$ is the unique shift-invariant measure reaching this supremum. In particular
we have
\begin{equation}\label{potentro}
h(\mu_\phi)=-\int \phi \dd\mu_\phi.
\end{equation}


\section{An exponential inequality and its general consequences}

\subsection{An exponential inequality}

Our main tool is an exponential inequality for fairly general observables.

Let $K:\Omega^n\to\rr$ be a function of $n$ variables and, for each $j=0,\ldots,n-1$,
let 
\begin{align*}
& \Lip_{j}(K) = \sup_{\x^{(0)}, \x^{(1)},\ldots,\x^{(n-1)}}\sup_{\y^{(j)}\neq\x^{(j)}}\\
& \frac{\left\lvert K\big(\x^{(0)},\ldots,\x^{(j-1)},\x^{(j)},\x^{(j+1)},\ldots,\x^{(n-1)}\big) -
K\big(\x^{(0)},\ldots,\x^{(j-1)},\y^{(j)},\x^{(j+1)},\ldots,\x^{(n-1)}\big)\right\rvert}
{d_{\theta}(\x^{(j)},\y^{(j)})}\cdot
\end{align*}
We shall say that  $K$ is a separately Lipschitz function of $n$ variables if
$$
 \Lip_{j}(K)<\infty,\; j=0,\ldots,n-1.
$$
We now present our main tool.

\begin{thm}[\cite{CMS}]\label{GaussIneqTh}
Let $\mu_{\phi}$ be a Gibbs measure. Then there exists a constant $D=D(\phi)>0$ such that,
for any integer $n\geq 1$ and for any separately Lipschitz function $K$ of $n$ variables, one has
\begin{equation}\label{GaussIneq}
\int e^{K(\x,\ldots,\sigma^{n-1}\x)}\ \dd\mu_{\phi}(\x)
 \leq e^{\int K(\y,\ldots,\sigma^{n-1}\y) \dd\mu_{\phi}(\y)} \ e^{D\sum_{i=0}^{n-1}\Lip_{i}^{2}(K)}.
\end{equation}
\end{thm}

Let us emphasize that the constant $D$ only depends on $\phi$. It depends
neither on $K$ nor on $n$.  \newline
The powerfulness of \eqref{GaussIneq} lies
in that it applies to {\em any} separately Lipschitz function of $n$ variables, regardless of its complicated
or implicit form. All we have to do is to estimate its Lipschitz constants.

The proof of Theorem \eqref{GaussIneqTh} easily follows from \cite{CMS} where it is done for a class of piecewise expanding maps of the interval (without assuming a Markov partition). In fact, the proof becomes simpler in our setting. It relies on the fact that the transfer operator associated to $\phi$
has a spectral gap when acting on Lipschitz functions. The point is to write an observable of
$n$ variables as a telescopic sum of observables depending only on one variable.

\subsection{General consequences}

We derive several consequences of inequality \eqref{GaussIneq}.\newline 
The first one is a bound for the probability of $K$ 
to deviate from its expectation.
\begin{coro}\label{ConcentrationIneq}
For every $t>0$, one has
\begin{equation}\label{pouic}
\mu_{\phi}\left\{ \x : K(\x,\ldots,\sigma^{n-1}\x)\geq
\int K(\y,\ldots,\sigma^{n-1}\y)\ \dd\mu_{\phi}(\y)+ t\right\} \leq e^{-\frac{t^{2}}{4D\sum_{i=0}^{n-1}\Lip_{i}^{2}(K)}}.
\end{equation}
\end{coro}

\begin{proof}
The proof is an immediate consequence of Markov inequality and \eqref{GaussIneq}: for every $\lambda>0$, 
the function $\lambda K$ is separately Lipschitz and
\begin{align*}
& \mu_{\phi}\left\{ \x : K(\x,\ldots,\sigma^{n-1}\x)\geq
\int K(\y,\ldots,\sigma^{n-1}\y)\ \dd\mu_{\phi}(\y)+ t\right\}  \\
& \qquad \leq   e^{-\lambda t}\ \int
e^{\lambda\big[K(\x,\ldots,\sigma^{n-1}\x)- \int K(\y,\ldots,\sigma^{n-1}\y) \dd\mu_{\phi}(\y)\big]}
\dd\mu_{\phi}(\x) \\
& \qquad \leq   e^{-\lambda t+\lambda^2D\sum_{i=0}^{n-1}\Lip_{i}^{2}(K)}.
\end{align*}
It remains to optimize over $\lambda$ to get the desired inequality.
\end{proof}

Of course we can apply \eqref{pouic} to $-K$ and get by a union bound that
$$
\mu_{\phi}\left\{ \x : \Big\vert K(\x,\ldots,\sigma^{n-1}\x)-
\int K(\y,\ldots,\sigma^{n-1}\y)\ \dd\mu_{\phi}(\y)\Big\vert\geq t\right\} \leq 2\
e^{-\frac{t^{2}}{4D\sum_{i=0}^{n-1}\Lip_{i}^{2}(K)}}
$$
for every $t>0$.

\bigskip

Another immediate consequence of \eqref{GaussIneq} is
a bound on the variance of $K$:
\begin{coro}\label{VarIneq}
One has 
$$
\int\Big(K(\x,\ldots, \sigma^{n-1}\x) - \int K(\y,\ldots,\sigma^{n-1}\y)\ \dd\mu_{\phi}(\y)\Big)^{2}\ \dd\mu_{\phi}(\x)
\leq 2D\sum_{i=0}^{n-1}\Lip_{i}^{2}(K).
$$
\end{coro}

\begin{proof}
We apply \eqref{GaussIneq} to $\lambda K$, with $\lambda\neq 0$, to get at once
\[
\frac{1}{\lambda^{2}}\left( \int e^{\lambda
\big[K(\x,\ldots,\sigma^{n-1}\x)- \int K(\y,\ldots,\sigma^{n-1}\y)\ \dd\mu_{\phi}(\y)\big]}\ \dd\mu_{\phi}(\x) -1\right) \leq \frac{1}{\lambda^{2}}\left(e^{\lambda^{2}D\sum_{i=0}^{n-1}\Lip_{i}^{2}(K)}-1\right).
\]
The result follows by Taylor expansion and letting $\lambda$ going to $0$. 
\end{proof}

The simplest, yet non-trivial, application of the above results is to ergodic sums, that
is to take $K_0(\x^{(0)}, \x^{(1)},\ldots,\x^{(n-1)})=f(\x^{(0)})+f( \x^{(1)})+\cdots +f(\x^{(n-1)})$
where $f:\Omega\to\rr$ is Lipschitz. 
A particular case of Corollary \ref{ConcentrationIneq} yields immediately
the following result, stated for later convenience.
\begin{coro}\label{phiphi}
Let $f:\Omega\to\mathbb{R}$ be a Lipschitz function. Then 
\begin{equation}\label{phiphiup}
\mu_{\phi}\left\{ \x : \frac{1}{n}\big(f(\x)+\cdots+f(\sigma^{n-1}\x)\big)-\int f \dd\mu_\phi\geq t\right\} \leq e^{-Bn t^{2}}
\end{equation}
for every $t>0$ and for every $n\geq 1$, where $B:=(4D |f|_\theta^2)^{-1}$.
\end{coro}
We can of course apply \eqref{phiphiup} to $-f$ to get
\begin{equation}\label{phiphidown}
\mu_{\phi}\left\{ \x : \frac{1}{n}\big(f(\x)+\cdots+f(\sigma^{n-1}\x)\big)-\int f \dd\mu_\phi\leq -t\right\} \leq e^{-Bn t^{2}}.
\end{equation}
By a union bound,  \eqref{phiphiup} and \eqref{phiphidown} yield at once
$$
\mu_{\phi}\left\{ \x : \left|\frac{1}{n}\big(f(\x)+\cdots+f(\sigma^{n-1}\x)\big)-\int f \dd\mu_\phi\right|\geq t\right\} \leq 2e^{-Bn t^{2}}.
$$
In words, the ergodic average of $f$ concentrates sharply around its $\mu_\phi$-average.
The bound is exponentially small in $n$ and when $t$
gets large, the probability of deviation is extremely small.

Let us close this section by a basic observation. Many estimators of interest are functions of
$n$ symbols, that is, functions of the form $\tilde{K}: A^{n}\to \rr$. A function $\tilde{K}: A^{n}\to \rr$
can be identified with a function $K:\Omega^n\to\rr$. When applying Theorem \ref{GaussIneqTh} and its corollaries
in this special case, $\Lip_j(K)$ has to be replaced by $\delta_j(\tilde{K})$, the oscillation at the $j$-th coordinate, where
\begin{align}
\label{lesdeltas}
&\delta_j(\tilde{K}) = \sup_{a_0,\ldots,a_{n-1}}
\sup_{a_j\neq b_j}\\
& \quad \quad\big| \tilde{K}(a_0,\ldots,a_{j-1},a_j, a_{j+1},\ldots,a_{n-1})-
\tilde{K}(a_0,\ldots,a_{j-1},b_j, a_{j+1},\ldots,a_{n-1})\big|.
\nonumber
\end{align}

\section{Bounds on entropy estimators}\label{clementine}

Throughout this section, $\phi\in \mathscr{F}_\theta$ and $\mu_\phi$ is its unique Gibbs measure.

\subsection{Plug-in estimator}

The plug-in estimator is based on the empirical frequency of a word $a_0^{k-1}$ in a `sample' $x_0,x_1,\ldots,x_{n-1}$:
$$
\E_{k}(a_{0}^{k-1};x_{0}^{n-1})=
\frac{1}{n}
\#\big\{0\leq j\leq n : \tilde{x}_j^{j+k-1}=a_{0}^{k-1} \big\},
$$
where $\tilde{\x}:=x_0^{n-1}x_0^{n-1}\cdots$ is the periodic point with period $n$ made from $x_0^{n-1}$.
This trick makes $\E_{k}(\cdot;x_{0}^{n-1})$ a locally shift-invariant probability measure on $A^k$.

For any ergodic measure $\nu$, there is a set of $\nu$-measure one of $\x$'s such that for every $k\geq 1$
$$
\lim_{n\to\infty}\E_{k}(a_{0}^{k-1};x_{0}^{n-1})=\nu([a_0^{k-1}]).
$$
The $k$-block empirical entropy is defined as
\[
\H_{k}(x_{0}^{n-1}):=H_{k}(\E_{k}(\cdot;x_{0}^{n-1})).
\]
It is clear that for $\nu$-almost every $\x$
$$
\lim_{k\to\infty} \lim_{n\to\infty} \frac{\H_{k}(x_{0}^{n-1})}{k}=h(\nu).
$$
As shown by Ornstein and Weiss (see \cite{shields}), we can in fact take a single limit by letting $k$ depend on $n$: 
if $k(n)\to\infty$ and $k(n)\leq \frac{1}{h(\nu)} \log n$ then
$$
\lim_{n\to\infty} \frac{\H_{k(n)}(x_{0}^{n-1})}{k(n)}=h(\nu)\quad\textup{for}\;\nu-\textup{almost-every}\;\x.
$$
Note that since $h(\nu)\leq \log|A|$ we can always take $k(n)\leq \frac{1}{\log|A|}\log n$.

\medskip

We can formulate our first result on fluctuations of the plug-in entropy estimator. We denote
by  $\ee$ the expectation and by $\Var$ the variance under $\mu_\phi$.

\begin{thm}\label{main1}
Let $D$ be the constant appearing in \eqref{GaussIneq}.
For every $\alpha\in(0,1)$, $t>0$ and $n\geq 2$ one has
$$
\mu_\phi\left\{\left\lvert \frac{\H_{k(n)}}{k(n)} -\ee\left(\frac{\H_{k(n)}}{k(n)}\right)\right\rvert\geq t \right\} \leq 
2\exp\left(-\frac{n^{1-\alpha}t^{2}}{16D(\log n)^2} \right)
$$
provided that $k(n)\leq \frac{\alpha}{2\log|A|}\log n$.\\
Moreover for every $n\geq 2$
$$
\Var\left(\frac{\H_{k(n)}}{k(n)}\right) \leq 8D \frac{(\log n)^2}{n^{1-\alpha}}\cdot
$$
\end{thm}

\begin{proof}
Given any integer $k\geq1$, consider the function $\tilde{K}:A^{n}\to \rr$ defined by
\[
\tilde{K}(s_{0},\ldots,s_{n-1}) = \widehat{H}_{k}(s_{0}^{n-1}).
\]
We estimate the $\delta_{j}(\tilde{K})$'s (see \eqref{lesdeltas} for the definition of $\delta_{j}(\cdot)$). We claim that
\[
\delta_{j}(\tilde{K}) \leq 2k\lvert A\rvert^{k}\frac{\log{n}}{n}, \quad\forall j=0,\ldots,n-1.
\]
Indeed, given any string $a_{0}^{k-1}$, the change of one symbol in $s_0^{n-1}$
can decrease $\mathcal{E}(a_{0}^{k-1};s_{0}^{n-1})$ by at most $k/n$. It is
possible that another string gets its frequency increased, and this can be 
at most by $k/n$. This is the worst case. We then use
the fact that for any pair of positive integers $l$ and $k$ such that $l+k\leq n$, one has
\[
\left\lvert \left(\frac{l}{n}\right)\log\left(\frac{l}{n}\right) - \left(\frac{l+k}{n}\right)\log\left(\frac{l+k}{n}\right)\right\rvert \leq \frac{k}{n}\log{n}.
\]
The claim follows by summing up this bound for all strings, which gives the factor $|A|^{k}$.\\
Finally, taking $k(n)\leq \frac{\alpha}{2\log\lvert{A}\rvert}\log{n}$, with $\alpha\in(0,1)$, and applying
Corollaries \ref{ConcentrationIneq} and \ref{VarIneq}, we get the desired
bounds.
\end{proof}

It is natural to seek for a concentration bound for the empirical entropy not
about its expectation, but about $h(\mu_\phi)$, the entropy of the Gibbs
measure. To have good control on this expectation, it turns out
that a better estimator is the conditional empirical entropy. To define
it, we need to recall a few definitions and facts.

For a shift-invariant measure $\nu$ and $k\geq 2$, let
\[
h_{k}(\nu) = H_k(\nu)-H_{k-1}(\nu)=-\sum_{a_{0}^{k-1}}\nu([a_{0}^{k-1}])
\log\frac{\nu([a_{0}^{k-1}])}{\nu([a_{0}^{k-2}])}\cdot
\]
It is well-known that $\lim_{k\to\infty} h_k(\nu)=h(\nu)$ (see for instance \cite{shields}).

The $k$-block conditional empirical entropy is 
$$
\hat{h}_k(x_0^{n-1})=\hat{h}_k(\E_{k}(\cdot;x_{0}^{n-1})).
$$
When $\nu$ is ergodic, one can prove \cite{shields} that, if
$k(n)\to\infty$ and $k(n)\leq \frac{(1-\epsilon)}{\log|A|}\log n$, for any
$\epsilon\in(0,1)$, then
$$
\lim_{n\to\infty} \hat{h}_{k(n)}(x_0^{n-1})=h(\nu),\quad\textup{for}\; \nu-\textup{almost every}\; \x.
$$

We have the following result.

\begin{thm}\label{main2}
Assume that $\theta<|A|^{-1}$.
There exist strictly positive constants $c,\gamma,\Gamma,\xi$ such that for every $t>0$ 
and for every $n$ large enough
$$
\mu_\phi\left\{\left\lvert \hat{h}_{k(n)} -h(\mu_\phi)\right\rvert\geq t + \frac{c}{n^{\gamma}}\right\} \leq 
2\exp\left(-\frac{\Gamma n^{\xi} t^{2}}{(\log n)^4} \right)
$$
provided that $k(n)< \frac{\log n}{2\log|A|}$.  
\end{thm}

\begin{rem}
From the proof we have
$\gamma=1/\big(1+\frac{\log |A|}{\log(\theta^{-1})}\big)$,
$\xi=1-2/\big(1+\frac{\log(\theta^{-1})}{\log |A|}\big)$ and
$\Gamma=(\log|A|)^2/16D$.
\end{rem}

\begin{proof}
By definition $\hat{h}_k=\H_k-\H_{k-1}$. If we let $\tilde{K}'(s_0,\ldots,s_{n-1})=\hat{h}_k(s_0^{n-1})$,
we estimate $\delta_j(\tilde{K}')$ by $2\delta_j(\tilde{K})$.\newline
We now estimate the expectation of $\hat{h}_{k(n)}$. We need the following lemma.
\begin{lemma}\label{decomp}
We have 
\begin{equation}\label{pif}
\hat{h}_{k(n)}(x_0^{n-1})= \frac{1}{n}\sum_{j=0}^{n-1}(-\phi(\sigma^j\x))+
\widehat{\Delta}_{k(n)}(x_0^{n-1})+\mathcal{O}(\theta^{k(n)})
\end{equation}
where 
\begin{equation}\label{ggg}
\big| \ee\big(\widehat{\Delta}_{k(n)}\big)\big|\leq \frac{M |A|^{k(n)}}{n},
\end{equation}
where $M>0$.
\end{lemma}
This lemma can be deduced from the proof of Theorem 2.1 in \cite{GGG}.
However, for the reader's convenience, we provide part of its proof in the appendix.

Now substract $h(\mu_\phi)$ and take the expectation on both sides of \eqref{pif},  to get, using
\eqref{potentro},
$$
\ee\big(\hat{h}_{k(n)}\big)-h(\mu_\phi)=\ee\big(\widehat{\Delta}_{k(n)}\big)+\mathcal{O}(\theta^{k(n)}).
$$
We now take $k(n)=q\log n/\log|A|$, where $0<q<1$ has to be determined. Choosing
$q=1/\big(1+\frac{\log\theta^{-1}}{\log|A|}\big)$ we easily get that
\begin{equation}\label{choco}
|\ee\big(\hat{h}_{k(n)}\big)-h(\mu_\phi)|\leq \frac{c}{n^\gamma},
\end{equation}
where $c>0$ is some constant and $\gamma=1/\big(1+\frac{\log |A|}{\log(\theta^{-1})}\big)$.

To end the proof, we apply Corollary \ref{ConcentrationIneq} and use \eqref{choco}. For the
exponent $\xi$ in the statement of the theorem be strictly positive, one must have $q<1/2$, which
is equivalent to the requirement that $\theta<|A|^{-1}$.
\end{proof}

\subsection{Hitting times}
Given $\x,\y\in\Omega$, let
$$
W_{n}(\x,\y) = \inf\{ j\geq 1 : y_{j}^{j+n-1}=x_{0}^{n-1}\}.
$$
Under suitable mixing conditions on the shift-invariant measure $\nu$, one can prove \cite{shields}
that
$$
\lim_{n\to\infty}\frac{1}{n}\log W_n(\x,\y)=h(\nu),\quad \textup{for}\;\nu\otimes\nu-\textup{almost every}\;(\x,\y).
$$
In particular, when $\nu$ is a Gibbs measure in the above sense, this result holds true \cite{ChU}.

We have the following concentration bounds for the hitting-time estimator.

\begin{thm}\label{thmWaiting}
There exist constants $C_1, C_2>0$ and $t_{0}>0$ such that, for every $n\geq 1$ and every $t>t_{0}$, 
\begin{equation}\label{positivePart}
(\mu_\phi\otimes\mu_\phi)\left\{(\x,\y):\frac{1}{n}\log W_{n}(\x,\y)>h(\mu_\phi)+t \right\} \leq C_1 e^{-C_2nt^{2}}
\end{equation}
and
\begin{equation}\label{negativePart}
(\mu_\phi\otimes\mu_\phi)\left\{(\x,\y):\frac{1}{n}\log W_{n}(\x,\y)< h(\mu_\phi)-t \right\} \leq C_1 e^{-C_2nt}.
\end{equation}
\end{thm}

\bigskip

Let us notice that the upper tail estimate behaves differently than the lower tail estimate
as a function of $t$. This asymmetric behavior also shows up in the
large deviation asymptotics \cite{ChU}.

Let us sketch the strategy to prove Theorem \ref{thmWaiting}.
We cannot apply directly our concentration inequality to the random variable $W_n$ for the
following basic reason. Given $\x$ and $\y$, the first time that one sees the first $n$ symbols of $\x$ in $\y$
is $W_n(\x,\y)$ and assume it is finite. If we make $\y'$ by changing one symbol in $\y$, we have a priori no control
on $W_n(\x,\y')$ which can be arbitrarily larger than $W_n(\x,\y)$ and even infinite. Of course,
this situation is not typical, but we are forced to use the worst case to apply our concentration inequality.
Roughly, we proceed as follows. We obviously have
$\log W_n=\log(W_n\mu_\phi([X_0^{n-1}]))-\log\mu_\phi([X_0^{n-1}])$. On the one hand, we use a sharp
approximation of the law of the random variables $W_n\mu_\phi([X_0^{n-1}])$
by an exponential law proved in \cite{Aba}. On the other hand, by the Gibbs property,
$\log\mu_\phi([x_0^{n-1}])\approx \phi(\x)+\cdots+\phi(\sigma^{n-1}\x)$ and we can use 
Corollary \ref{phiphi} for $f=\phi$.

\bigskip

\begin{proof}
We first prove \eqref{positivePart}. 
We  obviously have
\begin{align*}
& (\mu_\phi\otimes\mu_\phi)
\left\{(\x,\y):\frac{1}{n}\log W_{n}(\x,\y) > h(\mu_\phi) + t \right\} \\ 
& =(\mu_\phi\otimes\mu_\phi) \left\{(\x,\y):\frac{1}{n}\log W_{n}(\x,\y)+ \frac{1}{n}\log\mu_\phi([x_0^{n-1}]) 
- \frac{1}{n}\log\mu_\phi([x_{0}^{n-1}]) - h(\mu_\phi) > t \right\} \\
& \leq
(\mu_\phi\otimes\mu_\phi)\left\{(\x,\y):\frac{1}{n}\log\left[W_{n}(\x,\y)\mu_\phi([x_{0}^{n-1}])\right]>\frac{t}{2}\right\}\\
& \quad  +\mu_\phi\left\{\x:-\frac{1}{n}\log\mu_\phi([x_{0}^{n-1}]) - h(\mu_\phi) > \frac{t}{2}\right\}\\
& =:T_1+T_2. 
\end{align*}
We first derive an upper bound for $T_2$.\\
We use \eqref{BowenGibbs}, Corollary \ref{phiphi} applied to $f=-\phi$ and \eqref{potentro} to get
\begin{align*}
T_2 & \leq
\mu_\phi\left\{-\frac{1}{n}\big(\phi+ \cdots+\phi\circ\sigma^{n-1}\big)- h(\mu_\phi) > \frac{t}{2}-\frac{1}{n}\log C\right\}\\
& \leq  e^{-B nt^2}
\end{align*}
for every $t$ larger than $2\log C$.\newline
We now derive an upper bound for $T_1$.
To this end we apply the following result which we state as a lemma. It is an immediate consequence of
Theorem 1 in \cite{Aba}.
\begin{lemma}[\cite{Aba}]\label{Abadi}
Let 
$$
\tau_{[a_{0}^{n-1}]}(\y) := \inf\big\{ j\geq 1 : y_{j}^{j+n-1}= a_{0}^{n-1}\big\}\cdot
$$
There exist strictly positive constants $C,c, \lambda_{1}, \lambda_{2}$, with $\lambda_{1}<\lambda_{2}$, such that for every $n\in \nn$, every string $a_{0}^{n-1}$, there exists $\lambda(a_{0}^{n-1})\in[\lambda_{1},\lambda_{2}]$ such that 
$$
\left\lvert \mu_\phi\left\{ \y : \tau_{[a_{0}^{n-1}]}(\y) >\frac{u}{\lambda(a_{0}^{n-1})\mu_\phi([a_{0}^{n-1}])}\right\} - e^{-u}
\right\rvert\leq Ce^{-cu}
$$
for every $u>0$.
\end{lemma}

\medskip

By definition and using the previous lemma we get
\begin{align*}
T_1 &=\sum_{a_{0}^{n-1}}\mu_\phi([a_{0}^{n-1}])\
\mu_\phi\left\{ \y : \tau_{[a_{0}^{n-1}]}(\y)\mu_\phi([a_{0}^{n-1}])>e^{nt/2}\right\}\\
& \leq C'\ e^{-c' e^{nt/2}}
\end{align*}
for some $c',C'>0$.\newline
Since the bound for $T_1$ is (much) smaller than the bound for $T_2$, we can
bound $T_1+T_2$ by a constant times $e^{-B nt^2}$.
This yields \eqref{positivePart}.
 
We now turn to the proof of \eqref{negativePart}.
We have
\begin{align*}
& (\mu_\phi\otimes\mu_\phi)
\left\{(\x,\y):\frac{1}{n}\log W_{n}(\x,\y) < h(\mu_\phi) - t \right\} \\ 
& =(\mu_\phi\otimes\mu_\phi) \left\{(\x,\y):-\frac{1}{n}\log W_{n}(\x,\y)- \frac{1}{n}\log\mu_\phi([x_{0}^{n-1}]) 
+\frac{1}{n}\log\mu_\phi([x_{0}^{n-1}]) + h(\mu_\phi) > t \right\} \\
& \leq
(\mu_\phi\otimes\mu_\phi)\left\{(\x,\y):-\frac{1}{n}\log\left[W_{n}(\x,\y)\mu_\phi([x_{0}^{n-1}])\right]>\frac{t}{2}\right\}\\
& \quad +\mu_\phi\left\{\x:\frac{1}{n}\log\mu_\phi([x_{0}^{n-1}]) + h(\mu_\phi)  >\frac{t}{2}\right\}\\
& =T'_1+T'_2. 
\end{align*}
Proceeding as for $T_2$ (applying Corollary \ref{phiphi} to $f=\phi$) we obtain the upper bound
\begin{align*}
T'_2 & \leq
\mu_\phi\left\{\frac{1}{n}\big(\phi+ \cdots+\phi\circ\sigma^{n-1}\big)- \int\phi \dd\mu_\phi >\frac{t}{2}-\frac{1}{n}\log C\right\}\\
& \leq e^{-Bnt^2}
\end{align*}
for some $C">0$ and for every $t>2\log C$.

To bound $T'_1$ we use the following lemma (Lemma 9 in \cite{Aba}).
\begin{lemma} [\cite{Aba}]
For any $v>0$ and for any $a_0^{n-1}$ such that $v\mu_\phi([a_{0}^{n-1}])\leq 1/2$, one has
$$
\lambda_{1}\leq -\frac{\log\mu_\phi\left\{\tau_{[a_{0}^{n-1}]}>v\right\}}{v\mu_\phi([a_{0}^{n-1}])}\leq \lambda_{2},
$$
where $\lambda_{1}, \lambda_{2}$ are the constants appearing in Lemma \ref{Abadi}.
\end{lemma}
The previous lemma implies that 
$$
\mu_\phi\left\{\tau_{[a_{0}^{n-1}]}\mu_\phi([a_{0}^{n-1}])<v\right\}\leq 1-e^{-v\lambda_2}\leq \lambda_2 v
$$
provided that $v\mu_\phi([a_{0}^{n-1}])\leq 1/2$.
Taking $v=e^{-nt/2}$ it follows that
\begin{align*}
T'_1 &=\sum_{a_{0}^{n-1}}\mu_\phi([a_{0}^{n-1}])\
\mu_\phi\left\{ \y : \tau_{[a_{0}^{n-1}]}(\y)\mu_\phi([a_{0}^{n-1}])<e^{-nt/2}\right\}\\
& \leq \lambda_2\ e^{-nt/2}.
\end{align*}
This inequality holds if $e^{-nt/2}\mu_\phi([a_{0}^{n-1}])\leq 1/2$, which is the case for any $n\geq 1$ if $t\geq 2\log 2$.\\
Inequality \eqref{negativePart} follows from the bound we get for $T'_1+T'_2$.
But the bound for $T'_1$ is bigger than the one for $T'_2$, whence the result.\\
The proof of the theorem is complete.
\end{proof}

\appendix

\section{Proof of Lemma \ref{decomp}}

We start with the following identity:
\begin{align}
\nonumber
\hat{h}_k(x_0^{n-1}) & =-\sum_{a_0^{k-1}\in A^k} \E_k(a_0^{k-1};x_0^{n-1})
\log\frac{\E_k(a_0^{k-1};x_0^{n-1})}{\E_{k-1}(a_0^{k-2};x_0^{n-1})}\\
& =
\label{piupiu}
-\sum_{a_0^{k-1}\in A^k} \E_k(a_0^{k-1};x_0^{n-1}) \log\frac{\mu_\phi([a_0^{k-1}])}{\mu_\phi([a_1^{k-1}])}+
\widehat{\Delta}_{k}(x_0^{n-1}),
\end{align}
where 
$$
\widehat{\Delta}_{k}(x_0^{n-1}):=
$$
$$
-\sum_{a_0^{k-1}\in A^k} \E_k(a_0^{k-1};x_0^{n-1}) \log\frac{\E_k(a_0^{k-1};x_0^{n-1})}{\E_{k-1}(a_0^{k-2};x_0^{n-1})}
+
\sum_{a_0^{k-1}\in A^k} \E_k(a_0^{k-1};x_0^{n-1}) \log\frac{\mu_\phi([a_0^{k-1}])}{\mu_\phi([a_1^{k-1}])}=
$$
$$
-\sum_{a_0^{k-1}\in A^k} \E_k(a_0^{k-1};x_0^{n-1}) \log\frac{\E_k(a_0^{k-1};x_0^{n-1})}{\mu_\phi([a_0^{k-1}])}
+
\sum_{a_0^{k-1}\in A^k} \E_k(a_0^{k-1};x_0^{n-1}) \log\frac{\E_{k-1}(a_0^{k-2};x_0^{n-1})}{\mu_\phi([a_1^{k-1}])}=
$$
\begin{equation}\label{puce}
-H_k(\E_k(\cdot;x_0^{n-1})|\mu_\phi)+H_{k-1}(\E_{k-1}(\cdot;x_0^{n-1})|\mu_\phi),
\end{equation}
where 
$$
H_k(\eta|\mu_\phi)=\sum_{a_0^{k-1}\in A^k} \eta([a_0^{k-1}])\log\frac{\eta([a_0^{k-1}])}{\mu_\phi([a_0^{k-1}])}
$$
is the $k$-block relative entropy of $\eta$ with respect to $\mu_\phi$. The second term in \eqref{puce}
is equal to $H_{k-1}(\E_{k-1}(\cdot;x_0^{n-1})|\mu_\phi)$ because of the following two facts. First,
$\sum_{a_0\in A}\E_{k}(a_0^{k-1};x_0^{n-1})=\E_{k-1}(a_1^{k-1};x_0^{n-1})$. 
This is because $\E_{k}(;x_0^{n-1})$ is a locally shift-invariant probability measure on $A^k$.
Second,  $\sum_{a_{k-1}\in A}\E_{k}(a_0^{k-1};x_0^{n-1})=\E_{k-1}(a_0^{k-2};x_0^{n-1})$,
because the family $(\E_{k}(;x_0^{n-1}))_{k=1,2,\ldots}$ is consistent.

The quantity $|\widehat{\Delta}_{k}(x_0^{n-1})|$ is bounded above by $(M |A|^k)/n$ according to
\cite[formula (4.16)]{GGG}, where $M>0$ is a constant.

Now we deal with the first term in \eqref{piupiu}. We first introduce the function
$$
\phi_k(\y):=\log\frac{\mu_\phi([y_0^{k-1}])}{\mu_\phi([y_1^{k-1}])}
$$
which is a locally constant function on cylinders of length $k$.
It is easy to verify that $\|\phi-\phi_k\|_{\infty}\leq |\phi|_\theta \theta^k$ (this follows at once from \cite[Prop. 3.2 p. 37]{PP}).
We get that
$$
-\sum_{a_0^{k-1}\in A^k} \E_k(a_0^{k-1};x_0^{n-1}) \log\frac{\mu_\phi([a_0^{k-1}])}{\mu_\phi([a_1^{k-1}])}=
\frac{1}{n}\sum_{j=0}^{n-1} (-\phi(\sigma^j \x)) + \mathcal{O}(\theta^k).
$$
The proof of the lemma is complete.
\bibliographystyle{plain}
\bibliography{bib}

\end{document}